\documentclass{article}
\usepackage[utf8]{inputenc}
\usepackage{amsmath}
\usepackage{amssymb}
\usepackage{mathtools}

\newcommand{\defeq}{\vcentcolon=}
\newcommand{\eqdef}{=\vcentcolon}
\newcommand{\suchthat}{\ : \ }

\newcommand{\argmin}[1]{\arg\min_{#1} \ }
\newcommand{\optmin}[1]{\mathrm{opt}\min_{#1}}

\renewcommand{\star}[1]{#1^*}
\newcommand{\inv}[1]{#1^{-1}}

\renewcommand{\tilde}[1]{\widetilde{#1}}
\newcommand{\R}{\mathbb{R}}

\renewcommand{\epsilon}{\varepsilon}
\newcommand{\eps}{\varepsilon}

\newcommand{\grad}[2]{\nabla #1\left(#2\right)}
\newcommand{\hess}[2]{\nabla^2 #1\left(#2\right)}
\newcommand{\sumsq}[1]{\left\|#1\right\|^2}

\newcommand{\ltag}[1]{\label{#1}\tag{#1}}

\title{Critical Point Finding with Newton-MR by Analogy to Computing Square Roots}
\author{Charles Frye\\
Redwood Center for Theoretical Neuroscience\\
University of California, Berkeley}
\date{}

\begin{document}

\maketitle

\begin{abstract}

    Understanding of the behavior of algorithms for resolving the optimization problem
    (hereafter shortened to OP)
    of optimizing a differentiable loss function (\ref{OP1}), is enhanced by
    knowledge of the critical points of that loss function,
    i.e.~the points where the gradient is 0.
    Here, we describe a solution to the problem of finding critical points
    by proposing and solving three optimization problems:
    1) minimizing the norm of the gradient (\ref{OP2}),
    2) minimizing the difference between the pre-conditioned update direction and the gradient (\ref{OP3}),
    and
    3) minimizing the norm of the gradient along the update direction (\ref{OP4}).
    The result is a recently-introduced algorithm for optimizing invex functions, Newton-MR~\cite{roosta2018},
    which turns out to be highly effective at the problem of finding
    the critical points of the loss surfaces of neural networks~\cite{frye2019}.
    We precede this derivation with an analogous, but simpler, derivation
    of the nested-optimization algorithm for computing square roots
    by combining Heron's Method with Newton-Raphson division.
    
\end{abstract}

\section{Why Critical Points?}
The problem of approximate optimization of differentiable scalar functions is
\begin{equation}\ltag{OP1}
    \star{\theta} \defeq \argmin{\theta} L(\theta)
\end{equation}
for some ``loss function'' $L$.
This optimization is typically carried out by
\textit{local methods}, which make use of only evaluations of the function
and some finite number of its derivatives.
Prominent examples include gradient descent, Newton-Raphson, and (L-)BFGS~
\cite{boyd2004}.

What guarantees can be made about such algorithms?
Focusing on the simplest and most widely-used, stochastic gradient descent,
and narrowing to differentiable functions,
it can be shown~\cite{nesterov2004}
that the point of convergence, $\theta_\infty$, satisfies
the first-order local optimality criterion
\begin{equation}
    \grad{L}{\theta_\infty} = 0
\end{equation}

We call such a point a \textit{critical point} of the loss.
More concretely, it can be proven~\cite{nesterov2004}
that for any choice of $\eps > 0$,
there exists a finite number of iterations $T$
such that
\begin{equation}
    \sumsq{ \grad{L}{\theta_T} } < \eps
\end{equation}

At the most concrete, additional assumptions on
(Lipschitz-flavored) numerical properties of the function $L$
give asymptotic upper-bounding big-O rates, sometimes with matching lower bounds
(for details and proofs, see~\cite{boyd2004, nesterov2004}).

This is insufficient to prove that local methods ``work''
in the sense of solving the original optimization problem, \ref{OP1}.
Indeed, there are many examples where a gradient near or at $0$ is not a certificate of optimality.
For example, any point where the gradient vanishes but the Hessian is indefinite
(it has positive and negative eigenvalues) satisfies the first-order criterion but can be at
arbitrary height on the loss.

So we introduce a second-order optimality criterion:
\begin{equation}\label{eqn:second-order-opt}
    \lambda_{\min}\left(\hess{L}{\tilde{\theta}}\right) \geq 0
\end{equation}
where $\lambda_{\min}\left(M\right)$ of a matrix $M$ is its smallest eigenvalue.
Criterion~\ref{eqn:second-order-opt} has accompanying relaxation
\begin{equation}\label{eqn:second-order-opt-approx}
    \lambda_{\min}\left(\hess{L}{\tilde{\theta}}\right) > -\eps
\end{equation}
\noindent which eliminates from consideration points with (strongly) indefinite Hessians.
It can be shown~\cite{jin2017} that stochastic gradient descent
converges to points that satisfy~(\ref{eqn:second-order-opt-approx}).
This leaves as possible non-locally-optimal points of convergence only
those where the Hessian (approximately) vanishes in some directions
and third- or higher-order information is required to determine local optimality.
Information of order $N>3$ is hard or impossible to come by,
seeing as it involves generic $N$th-order tensors,
so we presume that we're living in a world where criterion~\ref{eqn:second-order-opt} is sufficient
for local optimality.

Then, the points which might cause stochastic gradient descent to fail
to optimize a function are those which
\begin{enumerate}
    \item satisfy the first-order local optimality criterion (have zero gradient)
    \item satisfy the second-order local optimality criterion (have a positive semi-definite Hessian)
    \item nevertheless have unacceptably high values of the loss.
\end{enumerate}
\noindent Such points are known as \textit{bad local minima}.
More generally, optimization algorithms may be attracted to all kinds of critical points,
and different kinds of critical points may interact heterogeneously with optimization algorithms,
e.g.~slowing down some but speeding up others.
A clear understanding of the nature of the critical points of the losses of a class of problems
can clarify empirical results on which algorithms perform best,
suggest the existence of superior algorithms,
and guide theory to produce better guarantees and tighter bounds.

In order to better understand the critical point-finding problem,
the problem of finding the roots, or zeroes,
of the gradient of the loss,
it is instructive to examine how the square root of a number is calculated
to high precision.

\section{The Square Root as an Optimization Problem}

Addition ($+$) and multiplication ($\times$) are simple operations, in the following sense:
given exact representations for two numbers, an exact representation for the results of $+$ and $\times$
applied to those two numbers can be obtained in finite time\footnote{
Specifically, addition is $O(n)$ and
multiplication is $O(n\log n)$~\cite{harvey2019},
courtesy of the Fast Fourier Transform and the convolution theorem.}.
The symbols $a + b$ and $a \times b$ represent the exact, finite output
of a concrete, finite-time algorithm.
That is, both operations define closed monoids over finite-length bit strings.

This is not true of division ($\div$), inverse ($\inv{}$), or square root ($\sqrt{}$).
In these cases, the operation is defined in terms of a promise regarding what happens
when the output of this operation is subjected to $\times$:
\begin{align}
    b = a \div c &\implies b \times c = a\label{eqn:div_promise}\\
    b = a^{-1} &\implies b \times a = 1\label{eqn:inv_promise} \\
    b = \sqrt{a} &\implies b \times b = a\label{eqn:sqrt_promise}
\end{align}
\noindent and for an exact representation of a number $a$,
the number that fulfills this promise might not have an exact representation,
as is famously the case for $\sqrt{2}$.
This makes algorithm design for these operations more complex than for $+$ and $\times$.

There are individual strategies for each,
but one general idea that turns out to be very powerful is
\textit{relaxation to an optimization problem}.
That is, we take the exact promises made above,
recognize them as statements of the optimality of the output for some criterion,
and then use that criterion to define an approximate promise,
true up to some tolerance $\eps$.

For $\sqrt{}$, we rearrange the promise~(\ref{eqn:sqrt_promise}), denoting its exact solution with a $\star{}$, into:
\begin{align}
    \sqrt{a} &\eqdef \star{b} \\
    a &= \star{b} \times \star{b}  \\
    0 &= \star{b} \times \star{b} - a  \\
    0 &= \sumsq{ \star{b} \times \star{b} - a }\label{eqn:sqrt_opt_criter}
\end{align}

\noindent and then recognize that, due to the non-negativity of the norm,
the last line (\ref{eqn:sqrt_opt_criter}) is a statement of optimality.
That is, the value $\star{b}$ that we are looking for is the
argument that minimizes the expression on the RHS:

\begin{equation}
   \star{b} = \argmin{b} \sumsq{ b \times b - a }
\end{equation}

Exactly minimizing this expression to arbitrary precision might be impossible,
so we consider instead the set, $\tilde{B}$ of all $\tilde{b}$s that make the expression
smaller than some criterion value $\eps$:

\begin{equation}\label{eqn:optmin_def}
    \tilde{B} \defeq \left\{\tilde{b} \suchthat \sumsq{ \tilde{b} \times \tilde{b} - a } \leq \eps \right\}
\end{equation}

This notation is unwieldy, so let's introduce a symbol,
$\optmin{}$, that means ``like an argmin, but only up to some $\eps > 0$",
abstracting away which exact member of the solution set (in (\ref{eqn:optmin_def}) above, $\tilde{B}$) is returned and which value of $\eps$ is chosen.
Because this is a solution, we include a $\star{}$; because it is inexact, we include a \  $\tilde{}$\ :

\begin{equation}
    \star{\tilde{b}} \defeq \optmin{b} \sumsq{ b \times b - a }
\end{equation}

The workhorse algorithm for solving these problems is \textit{Newton-Raphson}:
update the current estimate of a function's root by subtracting off
the value of that function divided by its derivative.
Using the fact that the derivative of $b \times b - a$ is $2b$, we obtain

\begin{align}
    b_{t + 1} &= b_t - \frac{b_t \times b_t - a}{2 b_t}\\
    &= \frac{b_t}{2} + \frac{a \div b_t}{2}
\end{align}

This method of calculating square roots is known as Heron's Method,
after Hero of Alexandria, who wrote it down in 60 AD,
some 16 centuries before Newton and Raphson would generalize it
to generic ($C_2$-smooth) root-finding problems
\cite{brown1999}.

The presence of a $\div$ may raise alarms\footnote{
But note that division by 2, shouldn't, since for binary numbers, it can be implemented by an $O(1)$ bitshift.},
since $\div$ was also among our troublesome operators.
It is in fact the case that $\div$ may \textit{also} be calculated with Newton-Raphson,
typically by computing $\inv{}$:
\begin{align}
    a \div b &= a \times \inv{b}\\
    \inv{b} &\eqdef \star{c} \\
    1 &= b \times \star{c}  \\
    0 &= c^{*-1} - b  \\
    0 &= \sumsq{ c^{*-1} - b }
\end{align}

Conveniently, this gives a Newton-Raphson update that is entirely
in terms of multiplications and additions:

\begin{align}
    c_{t + 1} &= c_t - \frac{\inv{c_t} - b }{-c_t^{-2}}\\
    &= c_t \times \left(2 - c_t \times b\right)
\end{align}

In the next section,
we define an analogous algorithm for finding critical points.
That is, we again try to solve a root-finding problem with Newton-Raphson,
but this introduces a division, which we reformulate as an optimization problem.

\section{Finding Critical Points Effectively}

Critical points are, in fact, roots of the system of equations
given by the gradient of the loss:

\begin{equation}
    \nabla L\left(\star{\theta}\right) = 0
\end{equation}

Inspired by our method for finding the root of the square,
we relax this root-finding problem to

\begin{equation}\ltag{OP2}
    \star{\tilde{\theta}} \defeq \optmin{\theta} \sumsq{ \grad{L}{\theta} }
\end{equation}

We could proceed by directly minimizing \ref{OP2},
introducing a surrogate loss $g(\theta) \defeq \sumsq{ \grad{L}{\theta} }$.
This method has been independently invented several times
over the preceding half-century
~\cite{mciver1972,angelani2000,pennington2017}.
But this surrogate loss is generically quadratically worse,
approximately, in condition number than is the original loss,
and is, ironically, subject to a ``bad local minimum'' property of its own~\cite{mciver1972,doye2002,frye2019}.

So we instead treat \ref{OP2} as a root-finding problem and apply Newton-Raphson.
Just as in the square root case, this introduces a division operation
(in disguise as a matrix inversion):

\begin{equation}
    \theta_{t+1} = \theta_t - \inv{\hess{L}{\theta_t}} \grad{L}{\theta_t}
\end{equation}

The first issue that arises is non-invertibility:
any Hessian with zero eigenvalues has no exact inverse.
This could be resolved by replacing the inverse, $\inv{}$,
with the Moore-Penrose pseudo-inverse~\cite{weisstein}.
But again, the more serious issue is that any inverse is defined in terms of an exact promise,
which we can only ever satisfy approximately.
The relevant promise here, introducing the symbol
$\star{\Delta \theta_t}$ for the ideal $\theta_{t+1} - \theta_t$,
is

\begin{equation}
    \hess{L}{\theta_t} \star{\Delta \theta_t} = - \grad{L}{\theta_t}
\end{equation}

\noindent with optimization-ready form

\begin{equation}\label{eqn:nr_promise}
    \sumsq{\hess{L}{\theta_t} \star{\Delta \theta_t} + \grad{L}{\theta_t}} = 0
\end{equation}

\noindent and therefore our \textit{third} optimization problem, \ref{OP3}, is:

\begin{equation} \ltag{OP3}
    \tilde{\star{\Delta \theta_t}} \defeq \optmin{\Delta \theta}
    \sumsq{\hess{L}{\theta_t} \Delta \theta + \grad{L}{\theta_t}}
\end{equation}

We can therefore solve \ref{OP2} by iteratively solving \ref{OP3}.
It would be very cute if \ref{OP3} were also solved using Newton-Raphson;
alas in fact it turns out to be preferable to use
a conjugate gradient-type method~\cite{choi2011}.
Note also that this problem only requires computation of the product
of the Hessian with a vector, rather than the explicit computation of the Hessian,
which can be made substantially faster~\cite{pearlmutter1994}.

But we are not done quite yet.
The optimality of the promise~(\ref{eqn:nr_promise})
can be derived for a quadratic function
(it essentially changes to basis in which the quadratic form is the identity matrix).
For non-quadratic functions, it is only approximately optimal,
in so far as the quadratic approximation to the function is close.

And so the $\Delta \theta_t$ given by solving \ref{OP3}
may not be the best choice.
It is, however, often a reasonable \textit{direction}
along which to search for good updates.
Therefore we redefine our update to be a scaled version
of the output of a solver for \ref{OP3}:

\begin{equation}
    \theta_{t+1} = \theta_t + \tilde{\star{\eta}} \tilde{\star{\Delta \theta_t}}
\end{equation}

\noindent and, as the notation suggests, define $\tilde{\star{\eta}}$ as a solution to
a \textit{fourth} optimization problem:

\begin{equation} \ltag{OP4}
    \tilde{\star{\eta}} \defeq \optmin{\eta \in \R^+} \sumsq{\grad{L}{\theta_t + \eta \tilde{\star{\Delta \theta_t}}}}
\end{equation}

\noindent Because it involves optimizing an $n$-dimensional system over a 1-dimensional parameter,
this is known as a \textit{line search problem}\footnote{Though the non-negativity restriction on $\eta$ makes it more accurately a \textit{ray} search problem.}, and we solve it using back-tracking line search~\cite{wolfe1969}.
That is, we start with a large initial guess $\eta_0$,
check whether the improvement for that step size is sufficient according to some criteria
(e.g.~the Wolfe conditions~\cite{wolfe1969}),
and, if not, reduce the step size by a multiplicative factor.

This particular combination,
applying \cite{choi2011} on \ref{OP3} for selecting an update direction,
then back-tracking line search to select a step size,
was recently proposed in~\cite{roosta2018}, under the name Newton-MR,
as a method for optimizing certain functions
whose critical points are all optima,
the \textit{invex} functions.
Newton-MR was recently
shown~\cite{frye2019}
to outperform other proposed
algorithms~\cite{dauphin2014,pennington2017}
on the problem of finding the critical points
of a linear neural network loss surface~\cite{baldi1989}.

\section{Conclusion}

In summary:
to better understand the optimization problem \ref{OP1},
we wish to find the critical points of the loss function.
We define the critical point finding problem as optimization problem \ref{OP2},
just as is done when computing square roots.
Also just as in the square root algorithm,
solving \ref{OP2} requires, in its inner loop,
another optimization problem, a form of division, be solved: \ref{OP3}.
The only additional complexity added by moving up to higher-dimensional problems,
like the loss surfaces of neural networks~\cite{frye2019},
is that the inner loop of \ref{OP2} gains another (non-nested) optimization problem,
\ref{OP4}, to select the step size.

\section*{Acknowledgements}
The author would like to thank
Neha Wadia, Nicholas Ryder, Ryan Zarcone, Dylan Paiton, and Andrew Ligeralde
for useful discussions.

\bibliography{bibliography}

\begin{thebibliography}{10}

\bibitem{roosta2018}
F.~Roosta, Y.~Liu, P.~Xu, and M.~W. Mahoney, ``Newton-{M}{R}: {N}ewton's method
  without smoothness or convexity,'' {\em arXiv preprint arXiv:1810.00303},
  2018.

\bibitem{frye2019}
C.~G. Frye, N.~S. Wadia, M.~R. DeWeese, and K.~E. Bouchard, ``Numerically
  recovering the critical points of a deep linear autoencoder,'' 2019.

\bibitem{boyd2004}
S.~Boyd and L.~Vandenberghe, {\em Convex Optimization}.
\newblock New York, NY, USA: Cambridge University Press, 2004.

\bibitem{nesterov2004}
Y.~Nesterov, {\em Introductory Lectures on Convex Optimization}.
\newblock Springer {US}, 2004.

\bibitem{jin2017}
C.~Jin, R.~Ge, P.~Netrapalli, S.~M. Kakade, and M.~I. Jordan, ``How to escape
  saddle points efficiently,'' 2017.

\bibitem{harvey2019}
D.~Harvey and J.~Van Der~Hoeven, ``{Integer multiplication in time O(n log
  n)}.'' Mar. 2019.

\bibitem{brown1999}
E.~Brown, ``Square roots from 1; 24, 51, 10 to dan shanks,'' {\em The College
  Mathematics Journal}, vol.~30, no.~2, pp.~82--95, 1999.

\bibitem{mciver1972}
J.~W. McIver and A.~Komornicki, ``Structure of transition states in organic
  reactions. general theory and an application to the cyclobutene-butadiene
  isomerization using a semiempirical molecular orbital method,'' {\em Journal
  of the American Chemical Society}, vol.~94, no.~8, pp.~2625--2633, 1972.

\bibitem{angelani2000}
L.~Angelani, R.~D. Leonardo, G.~Ruocco, A.~Scala, and F.~Sciortino, ``Saddles
  in the energy landscape probed by supercooled liquids,'' {\em Physical Review
  Letters}, vol.~85, no.~25, pp.~5356--5359, 2000.

\bibitem{pennington2017}
J.~Pennington and Y.~Bahri, ``Geometry of neural network loss surfaces via
  random matrix theory,'' in {\em International Conference on Learning
  Representations (ICLR)}, 2017.

\bibitem{doye2002}
J.~P.~K. Doye and D.~J. Wales, ``Saddle points and dynamics of
  {L}ennard-{J}ones clusters, solids, and supercooled liquids,'' {\em The
  Journal of Chemical Physics}, vol.~116, no.~9, pp.~3777--3788, 2002.

\bibitem{weisstein}
E.~W. Weisstein, ``Moore-penrose matrix inverse. {From MathWorld---A Wolfram
  Web Resource}.''
\newblock Visited on June 10, 2019.

\bibitem{choi2011}
S.-C.~T. Choi, C.~C. Paige, and M.~A. Saunders, ``{MINRES}-{QLP}: A {K}rylov
  subspace method for indefinite or singular symmetric systems,'' {\em {SIAM}
  Journal on Scientific Computing}, vol.~33, no.~4, pp.~1810--1836, 2011.

\bibitem{pearlmutter1994}
B.~A. Pearlmutter, ``Fast exact multiplication by the {H}essian,'' {\em Neural
  Computation}, vol.~6, pp.~147--160, 1994.

\bibitem{wolfe1969}
P.~Wolfe, ``Convergence conditions for ascent methods,'' {\em {SIAM} Review},
  vol.~11, pp.~226--235, Apr. 1969.

\bibitem{dauphin2014}
Y.~Dauphin, R.~Pascanu, {\c{C}}.~G{\"{u}}l{\c{c}}ehre, K.~Cho, S.~Ganguli, and
  Y.~Bengio, ``Identifying and attacking the saddle point problem in
  high-dimensional non-convex optimization,'' {\em CoRR}, vol.~abs/1406.2572,
  2014.

\bibitem{baldi1989}
P.~Baldi and K.~Hornik, ``Neural networks and principal component analysis:
  Learning from examples without local minima,'' {\em Neural Networks}, vol.~2,
  no.~1, pp.~53 -- 58, 1989.

\end{thebibliography}
\bibliographystyle{ieeetr}

\end{document}